\documentstyle{article}

\def \l {\lambda}

\def \E {{\bf E}}

\def \vep {\varepsilon}
\def \a {\alpha}

\def \z {\zeta}

\def \l {\lambda}

\def \L {\Lambda}

\def \t {\tilde}

\def \P {{\hbox{\,Pr}}}

\author{A.Khorunzhy\thanks{e-mail: khorunjy@math.uvsq.fr}\\
D\'epartement de Math\'ematiques \\
Universit\'e de Versailles -- Saint-Quentin\\ Versailles, FRANCE
}

\title{STOCHASTIC VERSION OF THE ERDOS-RENYI LIMIT 
THEOREM }

\date{}

\begin{document}


\maketitle

\begin{abstract}
We generalize the Erd\"os-R\'enyi limit theorem on the maximum
of partial sums of random variables to the case when the 
number of terms in these sums is randomly distributed.
Relations between this limit theorem and 
the spectral theory of  random
graphs and random matrices 
are discussed.

\end{abstract}

\vskip 0.2cm

\quad 
{\it Secondary:} 05C50, 05C80, 15A52

\section{Introduction} 

The Erd\"os-R\'enyi limit theorem
concerns the asymptotic behaviour of the random variables
$$
\eta(n,k) = \max_{i=1,\dots,n-k} S_i(k)/k,
\quad S_i(k) = 
\xi_i + \xi_{i+1} +\cdots + \xi_{i+k},
\eqno (1.1)
$$
where $\Xi = \{\xi_i\}_{i=1}^\infty$ is a family of independent
identically distributed (i.i.d.) random variables 
determined on the same probability space $\Omega$
and having zero mathematical expectation
$\E \xi=0$.
It is assumed  that the function
$$
\phi(\tau) = \E e^{\xi \tau}
\eqno (1.2)
$$
is determined  for $ \tau\in I_\xi$, where $I_\xi\subseteq {\bf R_+}=
(0,+\infty)$. 

In  \cite{ER} it is proved that  given $1<C<\infty$ there
exists, with probability 1, a  non-random limit 
$$
\lim_{n\to\infty} \eta(n,[C\log n]) = \a
\eqno (1.3)
$$
determined by relation
$$
\inf_{\tau \in I} \phi(\tau) e^{-\a\tau} = e^{-1/C}.
\eqno (1.4)
$$

In the particular case when $\xi_i$ are given as  Bernoulli random
variables
$$
\xi_i = \z_i = \cases{ 1, & with probability $1/2$, \cr
- 1, & with probability $1/2$,\cr}
\eqno (1.5)
$$
convergence (1.3) holds with $C= c\log_e 2$, where $\a(c)$ is
determined by relation 
$$
{1\over c} = 1- h\left( {1+\a \over 2}\right)
\eqno (1.6)
$$
with 
$$
h(t)  = - t\log_2 t - (1-t) \log_2(1-t), \quad 0<t<1.
$$ 
It is easy to see that  in this case $\a$ takes values between the mean
value
$0=\a(+\infty)$ and
the maximum $1=\a(1)$ of random variables $\z_i$. 
Obviously, one can also determine
the limit $\a(c)$  for the values $c\in (0,1)$; in this case it is
equal to $1$.

Further studies give more details about the convergence (1.3); in particular, 
the convergence in probability 
was proved in
sand the estimates with probability 1 were derived for the difference
${k\over \log k} [\eta(n,k) - \a(C)] $ 
\cite{DDL}. 
The Erd\"os-R\'enyi  theorem has found several applications
(see e.g. \cite{C,RV}) and
various generalizations of it have been considered
(random variables indexed by sets, \mbox{non-i.i.d.} random variables, 
random variables in Banach spaces and others).

One more version of this limit theorem is motivated
by the studies of spectra of  random matrices \cite{K}.
Namely, 
when regarding the weighted adjacency matrix of a random graph,
one observes that the spectral norm of such a matrix is bounded
from below by the maximum of the sums
$S_i(k)$ (1.1), where the number of 
terms $k$
is distributed at  random \cite{K}.
Then, in the limit of large dimension of such a sparse random matrix, 
one faces 
the problem that can be called 
the  stochastic version of the
Erd\"os-R\'enyi limit theorem.  
It is clear that in this direction one can find different
generalizations of the Erd\"os-R\'enyi theorem. 
In present paper we give the proof of the results
announced previously \cite{K1} in the form
maximally close to  (1.1).
We discuss other related settings at the end of the paper.

Let us complete this introduction
with expressions of gratitude to Profs. A. Rouault and E. Rio 
for the  interest to this work and valuable discussions.

\section{Main result and discussion}

Let us consider the family of i.i.d. random variables 
$\L = \{\l_i\}_{i=1}^\infty$ determined on the same probability space as $\Xi$, 
also
independent from
the family $\Xi$. These $\l_i$ 
take values in ${\bf N}$ according to the
law $\Pr\{ \l = k\} = q(k)$ such that $\E \l = p$. 
We assume  that the
function {\mbox{$
\psi_p(t) = \E e^{\l t} , \, \, t\in I_{\l} \subseteq {\bf R_+}
$ }}
exists,  
$$
\psi_p(t) =  e^{p\chi(t)(1+o(1))} \quad {\hbox{as}} \quad p\to\infty
\eqno (2.1)
$$
and $\chi(t)$ 
is analytic and satisfies conditions
$\chi(t) \ge 0$, $\chi(0) = 0$. 
It is easy to deduce from (2.1) that
$$
\P\{ \l\ge l \} \le \inf_{t\in I_\l} \psi_p(t) e^{- tl} =
e^{-p f(l/p)(1+o(1))}, \quad p\to\infty,
\eqno (2.2) 
$$ 
where 
$$
f(y) = \sup_{t\in I_\l} \left[yt - \chi(t)\right].
\eqno (2.3)
$$
We assume that 
$f(y)$ is the steep function, i.e. it takes value $+\infty$
when $y$ goes beyond the domain of definition of $f$.
Then by the G\"artner-Ellis theorem  (see e.g.
\cite{DZ})
$$
\P\{ \l \ge yp\} = e^{-p f(y)(1+o(1))}, \quad  y\ge 0, \quad  p\to\infty.
\eqno (2.4)
$$
Let us note that  $f(y)$ is non-negative, strictly convex  monotone
function. It attains its minimal value at the point $y' = (\E \l)/p = 1$.

\vskip 0.5cm 
{\bf Theorem 2.1} 
{\it Let us consider the sums
$$
 S_i(\l_i) = \xi_i + \xi_{i+1} + \dots + \xi_{i+\l_i},
\eqno (2.5)
$$
where $\{\xi_i\}$ are as in (1.1),
and determine
$$
\t \eta (n,p) = \max_{i=1,\dots, n} \eta_i(n,p),
\quad \eta_i (n,p) =  S_i(\l_i)/p.
$$
There exists with probability 1 a non-random limit
$$
\lim_{n\to\infty} \t \eta(n, C\log n) = \t \a,
\eqno (2.6)
$$
determined by the following relations:

{\bf i)} If 
$$
 D(\t\a/y)  = \max_{\tau \in I_\xi} \left[ {\t\a \tau \over y} - \log
\phi(\tau)
\right] ,
\eqno (2.7)
$$
then
$\t\a= \t\a(C)$ is determined by relation 
$$
\inf_{y\ge 0 }\left[f(y) +y D(\t\a/y)\right] = {1\over
C}\, 
\eqno (2.8)
$$
that generalizes (1.4);

{\bf ii)} In the case of Bernoulli random variables $\xi_i= \z_i$ (1.5), 
convergence (2.6) holds
with  $\t\a = \t\a(c)$ determined by relation (cf. (1.6))
$$
\inf_{y\in  (\t\a,+\infty)} \left\{ f(y) + y\left[ 
1- h\left({1\over 2} + {\t\a\over 2y}\right)\right]\right\} = {1\over c}\, ,
\eqno (2.9)
$$
where 
$ c = C/\log_e 2$.
}

\vskip 0.5cm

To compare this theorem with results of \cite{ER}, let us  consider 
first the case (ii) of Bernoulli 
 random
variables.
The next simplifying assumption
is that $\l_i$ have the Poisson distribution
with parameter $p$. This makes (2.5) close to the model arising in the studies
of sparse random matrices (see the end of this paper).
One can easily derive that in this case $\chi(t) = e^t-1, \, I_\l = {\bf R_+}$ 
and 
$$
f(y) = \cases {0, & if $ y \in (0,1) $, \cr
 y(\log y - 1) +1, & if $y\in [1,\infty) $.
\cr}
$$

The  function 
$$
g_a(y)= y\left[ 
1- h\left({1\over 2} + {a\over 2y}\right)\right]
\eqno (2.10)
$$
is positive and strictly decaying on $(a, +\infty)$;
the maximum is attained at $a$ and {\mbox{$g_a(a)=a$,}}
$g'_a(a) = -\infty$.
The solution of (2.9) exists for all $c\in (0,+\infty)$ 
and {\mbox{$\lim_{c\to\infty } \t\a(c) = 0$.}} This coincides with the value
of $\a(+\infty)=0$ given by (1.6). 

It is not hard to show that
{\mbox{$\t\a(c)> \a(c)$}}
for all finite values of $c$. 
Moreover, (2.9) implies that 
$ \t\a(c)$ infinitely increases as $c\to 0$.
This means that 
$$
\lim_{1\ll p \ll \log_2 n} \t \eta (n,p) = +\infty\, ,
\eqno (2.11)
$$
while the corresponding value of $\a(c), c\to 0$ remains equal to $1$.
This is an important  difference between the usual and stochastic
cases of the Erd\"os-R\'enyi limit theorem (see Section 4).

The reason for (2.11) is that in the limit $c\to 0$ the averaging in (2.5)
is not sufficient and $\tilde \eta(n,p)$ really 
searches for the maximum of variables $\eta_i$. 
This is provided by those variables that have  almost all $\z_i$  equal to $1$;
since 
one can see large deviations of the
the number of terms in $S_i(\l_i)$  with respect to $p$,
then one can obtain infinite values of $\t\eta(n,p)$ (2.11).

Thus, we conclude that 
the large fluctuations of $q(l)$ in the scale $p$ are responsible for (2.11).
This proposition is supported by the following observation.
Let us forget the Poisson distribution of  $\l$ and assume 
that there exists a finite   interval $Y\subset (0,\infty)$
such that  $ q([py])=o(e^{-p})$ for all $y \in \bar Y = {\bf R}\setminus Y$.
Then  we  determine
$f(y)$ as
$+\infty$ on $\bar Y$  and still consider (2.9).
In this case $\sup_{c}\t\a(c)$ is finite. Finally, we observe that 
if $f(y)$ is close to the
Dirac $\delta$-function {\mbox{$\delta(y-1)$,}} then $\t\a(c)$ is close to the 
values $\a(c)$ given by (1.5).

Summing up these arguments, we arrive at the conclusion that 
$\lim_{c\downarrow 0}\t\a(c)= \infty$  provided the fluctuations of
$\l_i$ around
$p$ are sufficiently large. 

It the general case of finite but unbounded random variables
$\xi_i$,
the  limit $\t\a(C)$ as $C\to 1$ can be infinite
already in the classical case of {\mbox{$\l_i \equiv k = [C\log n]$,}}.

\section{Proof of Theorem 2.1.}

As in \cite{ER}, we give the proof of the item (ii) 
concerning the Bernoulli random variables 
$\z_i $ and then describe  the changes needed to prove Theorem 2.1 
in the general case.

Let us show that for any positive $\vep$
$$
\P\{ \t \eta(n,c\log_2 n) \ge \t \a + \vep\} = O(n^{-\delta}), \quad n\to\infty,
\eqno (3.1)
$$
where $\delta>0$ depends only on $\vep$.
We start with  elementary inequality
$$
\P\{ \sup_{i=1,\dots,n} \eta_i(n,p) \ge x\} \le \sum_{i=1}^n  
\P\{ \eta_i(n,p) \ge x\} 
= n \P\{ \eta_1(n,p) \ge x\},
$$
where we used the fact that $\eta_i$ are identically
distributed. Observing that   
$\{\omega: \, \,\eta_1\ge x\} \subset 
 \{\omega: \, \, \l\ge px\} $,
we can write that
$$
\P\{\eta_1(n,p) \ge x\} = \sum_{l\ge px} q(l) \P\{ S(l)\ge px\}.
\eqno (3.2)
$$

Using the Stirling formula,
one can write that 
$$
{1\over 2^l} \sum_{(l+px)/2\le j\le l}
{{l}\choose {j}}
=
 \sum_{(l+px)/2\le j\le l}
{
 2^{-j \log_2 (j/l) - (l-j) \log_2(1-j/l)}\over \sqrt{2\pi j (1-j/l) -1}}(1+o(1)).
$$
Elementary computation shows that 
the last sum is estimated by its first term stimes a constant.
Then 
 we obtain  inequality
$$
\P\{ S(l)\ge px\} 
\le {U\over \sqrt l} \,\, 
2^{{- l\left[1- h\left({1\over 2} + {xp\over 2l}\right)\right]}}.
\eqno (3.3)
$$

If $x\ge  \t\a +\vep$,
then there exists $\delta>0$ that
$$
 f(y) +  y \left[ 
1- h\left({1\over 2} + {x \over 2y}\right)\right] = f(y) + g_x(y) > {1+\delta\over c}
\eqno (3.4)
$$
for all $y\ge 0$. 
It is clear that the minimal value of $f(y)$ is 
 $f(1)=0$. Since $f(y)$ is strictly convex and monotone,
it is continuous. Let us denote by $z$ the value such that 
$f(y) \le {\delta /( 2c)}$ for  $ 1\le y\le z$.
Then for all $0\le y \le z$ 
$$
cg_x(y)
 > 1+{\delta\over 2}.
\eqno (3.5)
$$
Using monotonicity of $g_x(\cdot)$,   
we derive from (3.3) inequality
$$
\P\{\eta \ge x\} \le 
\left(\sum_{ px\le l \le pz} + \sum_{pz\le l}\right)
{U\over \sqrt l} \,\, 
q(l)2^{{ l\left[h\left({1\over 2} + {xp\over 2l}\right)-1\right]}}
\le
$$
$$
U\sum_{xp\le l\le pz}  
2^{-p g_x(l/p)}
+
2^{-p g_x(z)(1+o(1))} \sum_{l\ge pz} q(l) .
$$
Taking into account that $p= c\log_2n$,
using (3.5) and combination of   (2.2) and (3.4), we obtain that 
$$
\P\{\eta_1(n,p) \ge x\} = O(n^{-1-\delta/4})
$$
because the number of terms in the first sum
is of the order $O(\log_2 n)$. Relation  (3.1) is proved.

To prove the almost sure estimate, we  follow the scheme of    \cite{ER}.
Let us 
consider the sequence of random variables 
$$
\t\eta_j\equiv
\t\eta(e^{(j+1)/C}-1,j).
$$ 
 Then (3.1) implies convergence of the series
$\sum_j\Pr\{\t\eta_j> \t\a\}$.
Now, taking into account that $\t\eta(n, C\log n) \le \t\eta_j$ for all $n$ such
that
$e^{j/C} \le n\le e^{(j+1)/C}-1$, we obtain relation 
$$
\Pr\{ \limsup_{n\to\infty} \t\eta(n,C\log n) \le \t\a \} = 1.
$$
This completes the estimate from above of $\lim \t \eta(n,c\log_2n)$ 
for the case of Bernoulli random variables.

In the general case one can  use 
inequality (see e.g. \cite{Cr})
$$
\P\{  S(l) \ge px\} = (2\pi l b_l)^{-1/2} e^{-lD(x/y)}\ , 
\eqno (3.6)
$$
where $0< b\le b_l \le B <\infty$, instead of (3.3). 
The remaining part of the proof repeats the arguments presented above.

\vskip 0.5cm

Now let us show that $\P\{\max_i \eta_i < \t\a -\vep'\}$ vanishes as
$n\to\infty$.
To do this, we take an integer $m$ and determine the subsets of $\Omega$  
$$
B_n(m) = \{ \omega\in \Omega: \sup_{i=1,\dots, n} \l_i < m\}
$$  
The next observation is that  the events 
$$
A_k(n,m) = \{\omega:  \eta_{km+1}(n,p)\le x \vert B_n(m) \}
$$ 
are jointly independent for all $0\le k\le  n(m)-1, \, n(m) = [n/m]$.
Thus, we can write that
$$
\P \{\sup_i \eta_i \le x\vert B_n(m)\} =  \prod_{ k=1}^{n(m)} \P \{A_k(n,m)\} 
$$
$$=
 \left( {\P\{ \eta_{1}\le x \cap B_n(m) \}\over \P(B_n(m)) }\right)
^{n(m)}, 
\eqno (3.7)
$$
where we denoted $\eta_{1} = \eta_{1}(n,p)$.
Regarding  elementary relations 
$$
\P\{ F\cap B_n(m) \} \le 1 - \P\{ \overline F\cap B_n(m)\}
$$
and
$$
\P\{ D\cap B_n(m)\} = 
\P\{ D\} - \P\{ D\cap \overline{B_n(m)}\} \ge 
\P\{ D\} - \P\{ \overline{ B_n(m)} \}
$$
with $F= \{ \omega: \eta_{1}\le x \}$ and $D= \bar F$, we can write that
$$
\P\{ \eta_{1}\le x \cap B_n(m) \}\le 1- \P\{ \eta_{1}>x\} + 
\P\{\overline{ B_n(m)} \}.
\eqno (3.8)
$$

Let us considers $\P\{ \eta_{1}>x\}$.
If $x < \t\a(c) - \vep$, then there exist $\delta'>0$ 
and  $z'>1$ that
$$
f(y) + y
\left[1- h\left({1\over 2} + {x\over 2y}\right)\right]
\ge {1-\delta'\over c}\, , \quad {\hbox{ for all }} \,   y \ge z_1
\eqno (3.9)
$$

The Stirling formula implies the following inequality inverse to 
(3.3)
$$
{1\over 2^l} \sum_{(l+px)/2\le j\le l}
{{l}\choose {j}}
\ge {u\over \sqrt l} \,\, 
2^{{ - l\left[1- h\left({1\over 2} + {xp\over 2l}\right)\right]}}.
\eqno (3.10)
$$
Using this estimate and remembering 
monotonicity of the function $g_x(\cdot)$
(2.10),  we derive from (3.2) relation
$$
\P \{ \eta_1> x\} \ge
n^{-c{ {z_1}\left[1- h\left({1\over
2} + {x\over 2z_1}\right)\right]}}
\sum_{ yÊ\ge pz_1} q(l).
$$
Now  (2.4) together with (3.9) imply  that
$$
\P \{ \eta_{1}>x\} = O(n^{-1+\delta'}).
\eqno (3.11)
$$ 
In the general case, one can use (3.6) 
instead of (3.10) and obtain (3.11).

Let us estimate $\P\{\overline{B_n(m)}\} \le  n \P\{\l_i\ge m\}$. We use
again (2.2) and observe that if $z''$ is such that 
$f(z'') \ge 3/C$, then  
$$
\P\{\overline{B_n(m)}\} = O(n^{-2}), \quad m = pz'', \quad n\to\infty.
\eqno (3.12)
$$
Now we can derive from (3.7), (3.8), (3.11) and (3.12) that
$$
\P \{\sup_i \eta_i \le x\vert B_n(pz'')\} \le \left( {1  - O(n^{-1+\delta'})\over 
1- O(n^{-2})}\right)^{n/(pz'')} = O(e^{-n^{\delta'/2}}).
$$
Finally, writing inequality
$$
\P \{\sup_i \eta_i \le x\}\le  
\P \{\sup_i \eta_i \le x\vert B_n(m)\} \P (B_n(m)) + 
\P\{\overline{ B_n(m)} \}
$$
 with $m=pz''$, 
we get  that  
$$
\P\{ \t \eta(n,p) < \t\a - \vep'\} = O(n^{-2}).
$$
Therefore $\P\{\liminf_{n\to\infty} \t\eta(n,C\log n) \ge \t\a\} = 1$.
This completes the proof of Theorem 2.1.
\hfill$\Box$

\section{Applications to random graphs and random matrices}

Let us consider the  adjacency matrix $A$ of   a simple  graph 
$\Gamma$ with the sets of vertices and edges denoted by $V$ and $E$, respectively.
If $\vert V\vert =N$ and the vertices are enumerated, 
then $A$ 
is an $N\times N$ real symmetric matrix with the entries 
$$
A^{(N)}_{ij} = \cases{1, & if the edge $e(i,j)\in E$,\cr
0, & if $e(i,j)\notin E$, \cr}, \quad i,j = 1,\dots, N.
\eqno (4.1)
$$
Often one calls the
set of eigenvalues of $A^{(N)}$ the spectrum of 
$\Gamma$ \cite{ZDS}. 

One of the  models of random graphs (see e.g. \cite{B}) is determined 
by the ensemble
$\{A^{(N,p)}\}$ of matrices whose entries $\{a_{ij}, \, i\le j\}$
are given as a family of jointly independent random variables
with distribution 
$$
a_{ij} = \cases { 1, & with probability $p/N$,\cr
0, & with probability $1-p/N$.\cr}
$$
Having a random graph $\Gamma^{(N,p)}$, one can ask about the asymptotic
behaviour of its  spectrum when $N\to\infty$, in particular, 
what  happens with the maximal (minimal) eigenvalue of $A^{(N,p)}$.
This question was addressed in \cite{K} in more general  setting than (4.1).

Namely, the random matrix ensemble $W_{ij}^{(N,p)} = a_{ij} w_{ij}$
has been studied, where $\{w_{ij}, \, i\le j \}$ are jointly independent 
random variables,
also independent from $\{a_{ij}\}$. It is assumed that the probability 
distribution of $w_{ij}$
has all  odd moments zero $m_{2k+1}=0$ and $m_{2k}\le k^{(1+\tau)k}$
with $\gamma\ge 0, k\ge 1$. 
Under these conditions, it was shown that
the spectral norm
of the matrix $\hat W^{(N,p)} = {1\over \sqrt p}W^{(N,p)}$ 
in the limit $N,p\to\infty$ converges  with probability 1 
to the limits
$$
\Vert \hat  W^{(N,p)}\Vert = \cases { 2v, & if $ p = O( (\log n)^{1+\gamma})$,\cr
+\infty, & if $ p = O( (\log n)^{1-\gamma}))$, \cr}  
\eqno (4.2)
$$
for any  $\gamma>\tau$. Here we denoted $v = \sqrt{\E w_{ij}^2},\,  i,j=1,\dots N$.

Slightly modifying computations of  \cite{K}, one can show that 
the same convergence (4.2) is valid for the spectral norm of 
$\Vert {1\over \sqrt p} A^{(N,p)}\Vert $ with $\gamma>0$.

To study the limit of 
$p= C\log N$, one has to carry out more accurate analysis than that
of \cite{K}. One of the possible results  can be obtained by using
the Theorem 2.1. Indeed,
one can write inequality
$$
\Vert \hat W^{(N,p)}\Vert^2 \ge \max_{i= 1,\dots,n} \Vert \hat W^{(N,p)}
e(i)\Vert^2 \equiv  \max_{i= 1,\dots,n} T_i(N,p),
$$
where $e(i)_j = \delta_{ij}$.  
Observing that
$$
T_i(N,p) = {1\over p} \sum_{j=1}^N a_{ij} w_{ij}^2 \ge 
{1\over p} \sum_{j\ge i}^N a_{ij} w_{ij}^2 \equiv \hat T_i(N,p),
$$
one faces the same problem as described in Theorem 2.1.
Indeed, 
 $p\hat T_i(n,p)$ 
is given by the sum of independent random variables and the 
number of terms is given by 
$\hat \l_i = \sum_{j= i}^N a_{ij}$ that approaches the Poisson
random variables $\l_i$ with parameters $ip/N\le p$, respectively. 
Thus  $p\hat T_i(n,p)$ 
resembles $S_i(\lambda_i)$ (2.5) with $\xi$ replaced by $\hat \xi_j = w_{ij}^2$. Let us
denote
$$
H(N,p) = \sup_{i=1,\dots, N}  T_i(N,p) \,\, {\hbox{and}}\,\,
\hat H(N,p) = \sup_{i=1,\dots, N} \hat T_i(N,p).
\eqno (4.3)
$$

So, the first difference between  (2.5) and (4.3) is that  
$$
\E \hat \xi_j =  v^2>0.
\eqno (4.4)
$$
However, it is easy to check that Theorem 2.1 remains valid 
in the case of (4.4).
Relations (2.6)-(2.8) do not change provided $\phi$ (1.2)
is replaced by $\hat \phi(\tau) = \E e^{\tau\hat \xi}$.
In this case $\hat \a (+\infty) = v^2$.

The following proposition is true that 
$$
\lim_{N\to\infty} H(N, C\log N) \le \hat  \a(C),
\eqno (4.5)
$$
where $\hat \a(C)$ is determined by (2.7) and (2.8) in terms of $\hat \phi(\tau)$.
We put inequality in (4.5) because the parameters of random variables 
$\hat \l_i$ are of the order $p$ provided $i \sim 1$ but decrease to zero 
when  $ i$ increases up to $N$. This is another difference between 
$\hat H(N,p)$ and $\tilde \eta (n,p)$ (2.5).

In this connection, it would be interesting to develop
the analogs of the Erd\"os-R\'enyi limit theorem
for maximums of $\hat T_i$ and  of $T_i$. 
It is natural to expect that $\lim H(N,C\log N) = \hat \a(C)$.
Of special interest is
the study of asymptotic behaviour of 
$\Vert \hat W^{(N,C\log N )}\Vert^2 $ also because 
in the limit $C\to\infty$ it is four times greater than that of $\hat \a(C)$.

Regarding the adjacency matrix $A^{(N,p})$, it is shown in \cite{KS} that
its 
maximal eigenvalue is closely  related
with the maximal degree $\Delta$ of a random graph.
Since the asymptotic behaviour of $\Delta$ is fairly well studied,
this gives an important source of information on the spectra of random graphs.
It could be interesting to find the limit of the spectral norm of $A^{(N,p)}$
in dependence on $C$, where  $p= C\log N, N\to\infty$.

The behaviour of sums  of the type (2.5) is interesting by itself 
in the following
aspect. Assume that the random variables $\l_i$ are such that $\E \l_i = p$
but the second moment of $\l_i$ does not exists. Then it is interesting
to know does the border $p\sim \log
n$ still remain to be the critical one for the maximums
$\tilde \eta(n,p)$.


\end{document}